\documentclass[12pt]{article}                              %% LaTeX2e

\usepackage{amsmath,amsthm,amsfonts,amscd,eucal}
\newlength{\dinwidth}
\newlength{\dinmargin}
\numberwithin{equation}{section}

\def\A{{\cal A}}
\def\B{{\cal B}}
\def\C{{\cal C}}
\def\D{{\cal D}}

\def\F{{\cal F}}
\def\G{{\cal G}}
\def\H{{\cal H}}
\def\I{{\cal I}}

\def\L{{\cal L}}
\def\M{{\cal M}}
\def\N{{\cal N}}

\def\Q{{\Omega}}
\def\R{{\cal R}}

\def\RR{{\mathbb R}}

\def\a{\alpha}
\def\b{\beta}
        \def\G{\Gamma}
        \def\D{\Delta}

\def\th{\vartheta}

       \def\L{\Lambda}

\def\e{\varepsilon}

\def\s{\sigma}

\def\f{\varphi}

\def\imply{\Rightarrow}

\def\ov{\overline}

\def\ad{\text{Ad}}

\newtheorem{Thm}{Theorem}[section]
\newtheorem{Cor}[Thm]{Corollary}
\newtheorem{Prop}[Thm]{Proposition}
\newtheorem{Lemma}[Thm]{Lemma}
\theoremstyle{definition}

\theoremstyle{remark}

\title{Conformal Nets,\\
Maximal Temperature\\
and Models from Free Probability}
\author{\textsc{Claudio D'Antoni$^{(1)}$, Roberto Longo$^{(1)}$, Florin
 Radulescu$^{(2)}$}\\{}\\
$^{(1)}$Dipartimento di Matematica, Universit\`a di Roma ``Tor
  Vergata''\\
  Via della Ricerca Scientifica, I-00133 Roma, Italy\\
  $^{(2)}$Department of Mathematics, University of Iowa\\
  Iowa City, IA 52246, U.S.A.}
 \footnotetext[1]{Supported in
 part by MURST and CNR-GNAFA.}

\date{September 19, 1998}

\begin{document}
\maketitle

\begin{abstract}We consider conformal nets on $S^1$ of von Neumann
algebras, acting on the full Fock space, arising in free probability.
These models are twisted local, but non-local.
We extend to the non-local case the general
analysis of the modular structure. The local algebras turn out to be
$III_1$-factors associated with free groups. We use our set up to show examples
exhibiting arbitrarily large maximal temperatures, but failing to satisfy the
split property, then clarifying the relation between the latter
property and the trace class
conditions on $e^{-\b L}$, where $L$ is the conformal Hamiltonian.
\end{abstract}

\vfill

\newpage
%%%%%%%%%%%%%%%%%%%%%%%%%%%%%%%%%%%%%%%%%%%%%%%%%%%%%%%%%%%%%
\markboth{C. D'Antoni, R. Longo and F. Radulescu}{Maximal Temperature}

\section{Introduction}
This note grew up as an attempt to combine ideas from Free
Probability  \cite{VDN} and Algebraic Quantum Field Theory
\cite{Haag}.

In this spirit, we construct conformal nets on $S^1$ of von Neumann
algebras acting on
the full Fock space, generalizing the single von Neumann algebra
construction in Free probability  \cite{S}. Such models
arise by second quantization with Boltzmann statistics
of the one-particle Hilbert spaces associated with derivatives of
the $U(1)$-current algebra. The local algebras are III$_1$ factors,
but not approximately finite-dimensional, indeed they have cores
isomorphic to $L(\mathbb F_{\infty})\otimes B(\H)$, where $L(\mathbb
F_{\infty})$
is the von Neumann algebra generated by the left regular
representation of the free group
$\mathbb F_{\infty}$ on infinitely many generators \cite{Rad,S}.

These nets are not local, but satisfy twisted locality, therefore
we are led to extend to the non-local case the
general analysis of conformal nets on $S^1$, in particular concerning the
geometric description of the modular structure, cf. \cite{BGL, FG}.

As a consequence, we clarify the relations between the trace
class property for the semigroup generated by the conformal Hamiltonian
and the split property, i.e. the statistical independence of the
observable von Neumann algebras associated to disjoint intervals
with positive distance.

More specifically, let $\A$ be a conformal net on $S^1$ with conformal
Hamiltonian $L$.
As is known, if $\A$ satisfies the trace class condition at all $\b>0$
 with $\textnormal{Tr}(e^{-\b L})\leq e^{a\b^{-r}}$ for some constants
 $a>0$ and $r>0$, then $\A$ is split \cite{BuWi}.

 Here we point out first that no requirement on the growth of
 $\textnormal{Tr}(e^{-\b L})$
 is necessary, namely the trace class condition $\textnormal{Tr}(e^{-\b
 L})<0, \ \forall\b>0$, alone implies the split property, and second
 that, on the other hand,
 the trace class condition at one fixed $\b>0$ is not sufficient.

The first statement is a rather direct consequence of results in
\cite{BDL}. Concerning the second one, we construct conformal nets on $S^1$
that violate the split property, although a trace class condition
for $e^{-\b L}$ is  satisfied with a maximal temperature $\b_0^{-1}$,
 namely
$$
\textnormal{Tr}(e^{-\b L})<\infty\Leftrightarrow \b > \b_0 ,
$$
where $\b_0$ may be arbitrarily small.

So far the Operator Algebras analysis of Conformal Quantum Field Theory
has been mostly restricted to the local case, namely to the
``observable algebra'' case. However non-local nets appear
naturally, for instance in the Fermionic case.

We thus take this opportunity to develop the general analysis somewhat
in more detail than necessary.
 Some of these results have the same proof as in the local
case and are mentioned for completeness.
Other results, as modular covariance, need however an adaptation.
\section{General properties of conformal nets on S$^1$}
By an interval $I$ of $S^1$ we shall always mean a {\it proper} interval,
namely $I$ and its complement $I'$ are assumed to have non-empty
interiors. We denote by $\I$ the set of intervals of $S^1$.

A {\it conformal precosheaf} (or {\it conformal net}\footnote{As $\I$ is not an
inductive family the
terminology ``precosheaf'' is more appropriate, but the name ``net''
has acquainted more familiarity.}) $\A$ of
 von~Neumann algebras on the intervals
of $S^1$ is a map
 $$
I\to\A(I)
 $$
from $\I$  to the set of  von~Neumann algebras on a Hilbert space $\H$
which verifies the following properties 1,2,3,4:
\begin{description}
	\item[$\textnormal{\textsc{1. Isotony}}:$] {\it If $I_1$, $I_2$ are
intervals and
 $I_1\subset I_2$, then
 $$
\A(I_1)\subset\A(I_2)\ .
 $$}
\end{description}
In the following the M\"obius group $PSL(2,\RR)$, hence its universal covering
$G$, acts as usual by
diffeomorphisms of $S^1$.
\begin{description}
	\item[\textnormal{\textsc{2. Conformal invariance}}:] {\it There is a
representation
$U$ of $G$ on $\H$ such that
 $$
U(g)\A(I)U(g)^*=\A(gI)\ ,
\qquad g\in G,\ I\in\I.
 $$}
\end{description}
\begin{description}
	\item[$\textnormal{\textsc{3. Positivity of the energy}}:$]
	{\it The generator of the
rotation subgroup $\th\to U(R(\th))$ is positive}, where $R(\th)$
denotes the (lifting to $G$ of the) rotation by an angle $\th$
(in the following we shall often write $U(\th)$ instead of
$U(R(\th))$).
\end{description}
 Let $I_0$ be the upper semi-circle. We identify as usual $I_0$ with the
positive real line $\RR_+$ via the Cayley transform and
 we consider the one parameter groups $\L_{I_0}(t)$ and
$T_{I_0}(t)$ of diffeomorphisms of $S^1$  that are conjugate by
the Cayley transform respectively to the dilations
$x\to e^{ t}x$ and translations $x\to x+t $ on $\mathbb R$.
 Moreover we consider the reflection of $S^1$ given by
 $r_{I_0}: z\to \bar z$
 where $\bar z$ is the complex conjugate of $z$.

For a general $I\in\I$  we choose $g\in G$ such that
$I=gI_0$ and set
$$
\L_I=g \L_{I_0} g^{-1},\quad
 r_I=g  r_{I_0} g^{-1},\quad
 T_I=g  T_{I_0} g^{-1}\ .
$$
 ($T_I$  is however well-defined only up to a rescaling of the parameter).

Recall the equivalence
between the positivity of the conformal Hamiltonian and the
positivity of the usual Hamiltonian energy, see e.g. Lemma B.5 in \cite{GuLo}.
\begin{description}
	\item[\textnormal{\textsc{4. Existence of the vacuum}}:]
	{\it There exists a unit $U$-invariant
vector $\Q$ (vacuum vector) which is cyclic for the von~Neumann algebra
$\vee_{I\in\I}\A(I)$ and separating for the von~Neumann
algebra $\cap_{I\in\I}\A(I)$}\footnote{If the precosheaf
$\A$ is local, i.e. the algebras associated with disjoint
intervals commute (see below), last property follows from the cyclicity
property.}. \end{description}
Notice that the {\it dual precosheaf}
$$
\hat\A(I)=\A(I')'
$$
 is a conformal net on $S^1$ with the same representation $U$ and,
since $\vee\hat\A(I)=\cap\A(I)$, $\Q$ is cyclic for
$\hat\A$. Of course
$$
\hat{\hat\A}=\hat\A\ .
$$

Let $r$ be an orientation reversing isometry of $S^1$ with $r^2=1$
(e.g. $r_{I_0}$). The action of $r$ on $PSL(2,\RR)$ by conjugation
lifts to an action $\s_r$ on $G$. We denote by $\G$ the
semidirect product of $G$ with $\mathbb Z_2$ via $\s_r$.
Since $\G$ is a covering of the group
generated by $PSL(2,\RR)$ and $r$, $\G$ acts on $S^1$. We call
(anti-)unitary a representation $U$ of $\G$ with operators on $\H$
such that $U(g)$ is unitary, resp. antiunitary, when $g$ is
orientation preserving, resp. orientation reversing.

The results in the following are known in the local case, but we
stress their independence of any local
commutativity assumption. The exposition follows in part \cite{GuLo}.
Where there are variations, we give a proof.

\begin{Thm}\label{genprop} Let $\A$ be a conformal precosheaf on $S^1$.
Then the following properties hold:
 \item[{$(i)$}]  \emph {Reeh-Schlieder theorem}: $\Q$
is cyclic and separating for each von Neumann algebra $\A(I)$,
$I\in\I$
 \item[{$(ii)$}] \emph {Modular covariance}:
 For any $I\in\I$ the modular group of
$\A(I)$ wuth respect to $\Q$ has the
geometric meaning corresponding to $\Lambda_I$, namely
$$
\Delta_{I}^{it}\A(I_0)\Delta_{I}^{-it} =\A(\L_I(-2\pi t)I_0), \
I_0\in\I, \ t\in\mathbb R,
$$
where $\Delta_I=\Delta_{\A(I)}$ denotes the modular operator
associated with $(\A(I),\Q)$;

the one-parameter group $z(t)=z_\A(t)$ of unitaries defined by
 $$
U(\L_I(-2\pi t))=\D_I^{it}z_\A(t)
 $$
commutes with $U(g)$,
$g\in G$, and belongs to the center of the gauge group\footnote{We
mean here the group of al unitaries $V$ on $\H$ such that $V\Q=\Q$
and $V\A(I)V^*=\A(I),\ I\subset\I$.}.

In particular the unitary, positive energy,
representation $U$ of $G$ is uniquely determined by
$\A$ by the formula
$$
U(\Lambda_I(2\pi t))=\Delta_I^{-i\frac{t}{2}}\Delta_{I'}^{i\frac{t}{2}} \ .
$$
 Moreover $U$ extends to an (anti-)unitary
representation of $\G$ determined by
\begin{equation}
		U(r_I)=J_I,\ I\in\I,\label{J}
\end{equation}
where  $J_I$ is the modular
conjugation associated with $(\A(I),\Q)$.
 If $g\in\G$ is orientation reversing, then
 \begin{equation}
 	\label{anti}
 U(g)\A(I)U(g)^*=\hat\A(gI)\ .
\end{equation}
$\A$ and $\hat\A$ have the same unitary representation of $G$
and
\begin{equation}
	\label{delta}
\Delta_{\hat\A(I)}^{it}=\Delta_{\A(I)}^{it}z_\A(2t)\ .
\end{equation}
 \item[{$(iii)$}] \emph{Additivity (and continuity)} :
 if $I$, $I_k$ are intervals,
and $I\subset\cup_k I_k$, then
 $$
\A(I)\subset\vee_k\A(I_k);
 $$
 if $\bar I$ denotes the closure of $I\in\I$ and $\bar I\supset \cap_k I_k$
 then, then
 $$
 \A(\bar I)\supset\cap_k\A(I_k) ,
 $$
 in particular $\A(I)=\A(\bar I)$.
 \item[{$(iv)$}] $U(2\pi)=\pm 1$
\end{Thm}
\begin{proof}
The proof of $(i)$ and $(iii)$ are as in the local case, see
\cite{FrJo}.

$(ii)$: First we observe that, for any $g\in G$, $I\in\I$,
\begin{align*}
\D_I^{it}U(g)\D_I^{-it}&=U(\L_I(-2\pi t)g\L_I(2\pi t))\\
J_IU(g)J_I&=U(r_I g r_I)\ .
\end{align*}
 These relations
hold for $g=T_I(s)$ and for  $g=T_{I'}(s)$ by Borchers' theorem \cite{Bo},
hence for any
$g$ because $T_I$ and $T_{I'}$ generate $G$.
 Then, if $I_1\in\I$, we may find $g\in G$ such that
$I_1=gI$, therefore
\begin{align*}
\D_I^{it}\A(I_1)\D_I^{-it}
&=\D_I^{it}\A(gI)\D_I^{-it}\\
&=\D_I^{it}U(g)\D_I^{-it}
  \A(I)\D_I^{it}U(g^{-1})\D_I^{-it}\\
&=U(\L_I(-2\pi t)g\L_I(2\pi t))\A(I)
  U(\L_I(-2\pi t)g^{-1}\L_I(2\pi t))\\
&=\A(\L_I(-2\pi t)g\L_I(2\pi t)I)=\A(\L_I(-2\pi t)I_1)\ .
\end{align*}
 For any $t\in\RR$, $z(t)=\D_{I}^{it}U(\L_{I}(2\pi t))$ is thus an
automorphism of any local algebra, hence commutes with
$\D_{I}^{is}$, which implies that $t\to z(t)$ is a one parameter
group, commutes with $U(g)$, $g\in\G$, and is independent of
$I$. Due to the same reason $z(t)$ commutes with $J_{I}$.
 Then, setting $U(r_I)=J_I$, $g\in\G\to U(g)$ is an
(anti-)unitary representation which commutes with $z(t)$.

 Finally, if $g$ is orientation reversing, then
$gr_I$ is orientation preserving, therefore
\begin{multline*}
\ad U(g)\A(I)=\ad U(gr_I r_I)\A(I)=\ad U(gr_I)\ad J_I\A(I)\\
=\ad U(gr_I)\hat\A(I')=\hat\A(gr_I I')=\hat\A(gI)\ .
\end{multline*}
By uniqueness, the representation $U$ of $\A$ is also the
representation associated with $\hat\A$.

To show (\ref{delta}), notice first that $z_\A(t)$ commutes
with $J_I$, as it implements an automorphism of $\A(I)$,
hence with  $U(g),\ g\in\Gamma$ because
of its definition in (\ref{J}). Let then $g\in\Gamma$
be orientation reversing with
$gI =I$, thus $g\L_I(t)g^{-1}=\L_I(-t)$. We have
\begin{multline*}
\Delta_{\hat\A(I)}^{it}=U(g)\Delta_{\A(I)}^{-it}U(g)^*
=U(g)U(\L(2\pi t))z_\A(t)U(g)^*\\
=U(\L(-2\pi t))z_\A(t)=\Delta_{\A(I)}^{it}z_\A(2t).
\end{multline*}
The above formula also entails that $z(t)$ commutes with every
unitary $V$  in the gauge group, as such a $V$ commutes both with
$\Delta_{\A(I)}^{it}$ and $\Delta_{\hat\A(I)}^{it}$ by the
modular theory.

$(iv)$: Since $U(2\pi)$ is an automorphism of $\A(I)$, it
commutes with the associated modular antiunitary $J_I$.
Then,
 $$
U(2\pi)=J_IU(2\pi)J_I
=U(r_{I}R(2\pi)r_{I})=U(-2\pi)\ .
 $$
 \end{proof}
 \emph{Remark.} The group $G$ acts on on the $n$-covering and on the
 universal covering of $S^1$ (homeomorphic to $\mathbb R$) and one
 may consider more general precosheaves on these spaces. There are
 direct extensions of the above results in these cases, in particular
 concerning modular covariance.
 We omit this generalization for simplicity.

 \smallskip
  We shall say that $\A$ satisfies {\it twisted locality} if there
 exists a unitary $Z$, commuting with the unitary representation $U$
 and with $z_\A(t)$, such that $Z\Q=\Q$ and
  $$
  Z\A(I')Z^*\subset \A(I)'
  $$
  for all intervals $I$.
 \begin{Prop} Let $\A$ satisfy twisted locality. Then
 \item[(i)]  \emph{Twisted duality} holds:
 $$
   Z\A(I')Z^* =\A(I)', \quad I\in\I\ .
  $$
i.e. $\hat\A(I)=Z\A(I)Z^*, \ I\in\I$,
\item[(ii)] The \emph{Bisognano-Wichmann property} holds:
$$
\Delta_I^{it}=U(\L_I(-2\pi t)),
$$
namely $z(t)=1$.
  \end{Prop}
\begin{proof} $(i)$: By twisted locality $Z\A(I)Z^*\subset\hat\A(I) $,
moreover $\Q$ is cyclic and separating for both $\hat\A(I)$
and $Z\A(I)Z^*$. Now the modular group Ad$\Delta_{\hat\A(I)}^{it}$ of
$\hat\A(I)$
leaves $Z\A(I)Z^*$ globally invariant by equation (\ref{delta}) and the
commutativity
between $Z$ and $U$,$z(t)$; hence $\hat\A(I)=Z\A(I)Z^*$ by Takesaki's theorem.

$(ii)$: We have $z(t)=\Delta_I^{it}U(\L_I(2\pi t))$ independently of the
interval $I$, hence
$z(t)=\Delta_{I'}^{it}U(\L_{I'}(2\pi t))=\Delta_{I}^{-it}U(\L_{I}(-2\pi
t))=z(-t)$,
thus $z(t)=1$, where we have used the twisted duality property
to entail that
$\Delta_{I'}^{-it}=Z\Delta_{I}^{it}Z^*=\Delta_{I}^{it}$, as
$Z$ commutes with $\Delta_{I'}^{it}$.
\end{proof}
  Since $U(2\pi)$ is an involutive automorphism of any local algebra and has
square $1$, we may define the Fermi and Bose part of $\A(I)$ as
 $$
\A_{\pm}(I):=\{x\in\A(I):U(2\pi)xU(2\pi)=\pm x\}\ .
 $$
 As is  known and easy to check, normal commutation relations
are equivalent to twisted locality with the unitary $Z$ given by
\begin{equation}
Z=\frac{I+iU(2\pi)}{1+i}\ .\label{norm}
\end{equation}
 In the following proposition we show that a weak form of twisted
locality holds, i.e. the vacuum expectations of the commutators
vanish.
\begin{Prop} Let $\A$ a conformal precosheaf on $S^1$.
Then \emph{weak twisted locality} holds i.e., for any $I\in\I$,
 $$
([x,ZyZ^*]\Q,\Q)=0\ ,\qquad x=x^*\in\A(I),\ y=y^*\in\A(I'),
 $$
 where $Z$ is given by formula (\ref{norm}). In particular
 \emph{weak locality} (i.e. $([x,y]\Omega,\Omega)=0$)
 is equivalent to $U(2\pi)=1$
 \end{Prop}
 \begin{proof}
 By conformal invariance, it is sufficient to prove the weak
twisted locality for the upper semicircle $I_0$. This amounts to show that
when $x$ is a selfadjoint element in $\A({I_0})$ and $y$ is a
selfadjoint element in $\A(I'_0)$ then $(x\Q,Zy\Q)$ is real.

Since $Z$ commutes with $U(g)$ when $g$ is orientation preserving
and $JZJ=Z^*$, and making use of the commutation relations
following from Theorem \ref{genprop}, a straightforward computation
shows that, if $S_0$ is the Tomita operator for $(\A(I_0),\Q)$,
then
 $ZU(\pi)S_0 U(-\pi)=S_0^*Z$. Therefore we get
\begin{align}
(x\Q,Zy\Q)&=(x\Q,ZU(\pi)S_0U(-\pi)y\Q)=(x\Q,S_0^*Zy\Q)\\
&=(Zy\Q,x\Q)=\ov{(x\Q,Zy\Q)}\ .
\end{align}
 Finally let us assume weak locality and suppose by contradiction
that $U(2\pi)\ne1$. Then there exists $I\in\I$ and a non-zero
selfadjoint $x\in\A(I)$ such
that $U(2\pi)x\Q=-x$, hence, for any $y=y^*\in\A(I')$, $(x\Q,y\Q)$ is
real by weak locality and $(Zx\Q,y\Q)=i(x\Q,y\Q)$ is real too,
namely $(x\Q,y\Q)=0$ for any $y\in\A(I')$, and this implies the
thesis by the Reeh-Schlieder property.
 \end{proof}
We shall say that $\A$ is \emph{irreducible} if
the von~Neumann algebra $\vee\A(I)$
generated by all local algebras coincides with $\B(\H)$.
The irreducibility property is indeed equivalent to several other
requirement.
 \begin{Prop}\label{irr} Assume $z_{\A}(t)=1$. The following are equivalent:
 \item[{$(i)$}] $\mathbb C\Q$ are the only $U(G)$ invariant vectors.
\item[{$(ii)$}] The algebras $\A(I)$, $I\in\I$, are
factors. In particular they are type III$_1$ factors, provided
$\A\neq\mathbb C$.
\item[{$(iii)$}] The net $\A$ is irreducible.
\item[{$(iv)$}] The dual net $\hat\A$ is irreducible, i.e. the algebra
$\cap\A(I)$ given by the intersection
of all local algebras coincides with $\mathbb C$.
 \end{Prop}
\begin{proof}
 The proof is similar to the one given in \cite{GuLo} in the local
 case; one just notices that property $(i)$ is the same for $\A$
 and $\hat\A$, thus $(iii)$ and $(iv)$ are equivalent.
\end{proof}
 Irreducibility is also
equivalent to $\Q$ being unique invariant for any of the unitary subgroups
corresponding to $T_I$, $\Lambda_I$ or $R$, see
Lemma B.2 of the appendix \cite{GuLo}.

In the next Corollary the  assumption of compactness
$\{\textnormal{Ad}z_\A(t),\ t\in\mathbb R\}^-$ is satisfied in
particular if $\A$ is distal split \cite{DL} (see below) or, of
course, if $\A$ is twisted local.
\begin{Cor}\label{factor} If $\A\neq\mathbb C$ is irreducible and
$\{\textnormal{Ad}z_\A(t),\ t\in\mathbb R\}^-$ is compact,
then $\A(I),\ I\in\I$, is a type III factor. If $z(t)=1$, then
$\A(I)$ is of type III$_1$.
\end{Cor}
\begin{proof}
By assumption the closure (in the gauge automorphism group) ${\cal G}_0$ of
$\{\textnormal{Ad}z_\A(t)$, $t\in\mathbb R\}$ is a compact and abelian.
Denote by $\A^z$ the fixed-point net under the
action of ${\cal G}_0$. If $\A^z(I)=\mathbb C$, namely ${\cal G}_0$ acts
ergodically, then $\A(I)$ must be abelian, hence $\A$ is local and
$z_\A(t)=1$ (Prop. \ref{local}), but this implies that $\A^z=\A$ is trivial.

So we may assume that $\A^z$ is non-trivial,  hence $\A^z(I)$ is a type
III$_1$ factor by the uniqueness of the vacuum and Proposition
\ref{irr}. By a similar reasoning the relative commutant ${\cal R}(I)=
\A^z(I)'\cap\A(I)$ is abelian as ${\cal G}_0$ acts ergodically on it,
hence $\cal R$ has to be a constant local precosheaf. Then
$U(\L_I(-2\pi t))$ and $z_\A(t)$ have the same restriction to
the Hilbert space of $\cal R$, thus ${\cal G}_0$ would act trivially
on $\cal R$ (because  by lemma \ref{reg} $U(\L_I(-2\pi t))$
has no eigenvalue other than 1, while $z(t)$ has pure point spectrum)
 and ${\cal R}(I)=\mathbb C$. It follows that $\A(I)$ is a
factor that must be of type III because it has a normal conditional
expectation onto $\A^z(I)$.
\end{proof}
We include the following corollary, a variant of  a result in \cite{EvBo}.
 \begin{Cor}
Let $\A$ be an irreducible conformal net with $\{z_\A(t),\ t\in\mathbb
R\}^-$ compact and let
 $\A_0$ be its restriction to $\mathbb R = S^1\backslash \{-1\}$.
 Let $\B_0\subset\A_0$ be a subnet on $\mathbb R$,
  with the same translations and
dilations of $\A_0$, Then $\B_0(\mathbb R^+)'\cap\A_0(\mathbb R^+)$
is either trivial or a type III factor. If moreover $\B_0\subset\A_0$
has finite index, i.e.
$[\B_0(\mathbb R^+):\A_0(\mathbb R^+)]<\infty$, then the first case
occurs: $\B_0(\mathbb R^+)'\cap\A_0(\mathbb R^+)=\mathbb C$.
 \end{Cor}
 \begin{proof}
We assume that $\A$ is local; the general case follows by similar
arguments. As the translations
 and dilations of $\A_0$ restrict to $\B_0$,
it follows that $\B_0$ extends to a conformal precosheaf on $S^1$ \cite{GLW}.
Moreover
the uniqueness of the vacuum holds
for $\A_0$ hence for $\B_0$ too and therefore $\B(\mathbb R^+)$
is type III unless $\A$ is trivial (Prop. \ref{factor}).
By the same reason
$\B_0(\mathbb R^+)'\cap\A_0(\mathbb R^+)$ is either trivial or a type III
factor,
but the last possibility cannot occur in the finite-index case because
$\A_0(\mathbb R^+)'\cap\A(\mathbb R^+)$ is then
finite-dimensional.
 \end{proof}
\smallskip
We shall say that $\A$ is \emph{local} if whenever $I_1$, $I_2$ are disjoint
intervals the two algebras $\A(I_1)$ and $\A(I_2)$ commute.

\smallskip\noindent
\emph{Remark.} The ``observable net'' $I\to \C(I)=\A(I)\cap\hat\A(I)$ is local
 and, by modular covariance
 and Takesaki's theorem, for each $I$ there exists a normal,
vacuum-preserving, conditional
expectation $\e_I$ from $\A(I)$ (or from $\hat\A(I)$ onto $\C(I)$)
and $\e_{\tilde I}|_{\A(I)}=\e_I$ if $I\subset\tilde I$.
\begin{Prop}\label{local} Let $\A$ be a conformal precosheaf on $S^1$.
The following are equivalent:
\begin{itemize}
 \item[{$(i)$}] $\A$ is local.
\item[{$(ii)$}] Haag duality holds, i.e.
$\A(I')=\A(I)'\ ,\quad I\in\I.$
\item[{$(iii)$}] The vacuum is cyclic for the algebra
$\A(I)\cap\hat\A(I)$ for some (hence for any) $I\in\I$.
\item[{$(iv)$}] The algebras $\A(I)$ and $\hat\A(I)$ coincide
for some (hence for any) $I\in\I$.
\end{itemize}
If these properties hold then $U$ is indeed a representation of
$PSL(2,\mathbb R)$, i.e. $U(2\pi)=1$.
 \end{Prop}
 \begin{proof} The relations $(ii)\imply(i)$, $(iv)\imply(iii)$,
$(iv)\Leftrightarrow(ii)$ are obvious, and the implication
$(i)\imply(ii)$ is proven e.g. in \cite{BGL}. If $(iii)$
holds, then $\A(I)\cap\hat\A(I)$ is a subalgebra of $\A(I)$ for which
$\Q$ is cyclic. Since $\Delta_{\hat\A(I)}^{it}=\Delta_{\A(I)}z(2t)$,
we have that $\A(I)\cap\hat\A(I)$ is stable under the modular
group of $\hat\A(I)$ with respect to $\Q$, hence coincides with $\hat\A(I)$.
By the same argument it coincides with $\hat\A(I)$, too. Finally if the net is
local it is in particular weakly local, and this, according to
Proposition~1.2, implies the thesis.
\end{proof}
\subsection{Wiesbrock's theorem in the non-local case}
 In order to clarify the general structure, we give here the necessary
 modifications in order to characterize
conformal nets on $S^1$, also in the non-local case, by a simple
extension of Wiesbrock's theorem with \emph{half-sided modular inclusions}
(hsm).  This is however independent of the rest of the paper.

Recall that $(\R\subset{\cal S},\Q)$ is a $\pm$hsm if
$\R\subset{\cal S}$ are von
Neumann algebra, $\Q$ is a cyclic and separating vector for both $\R$
and ${\cal S}$ and
Ad$\Delta_{\cal S}^{it}{\cal S}\subset{\cal S},\ \pm t>0$,
where $\Delta_{\cal S}$
is  the modular operator associated with $({\cal S},\Q)$.

Let $\N$ and $\M$ be von Neumann algebras on a Hilbert space $\H$ and
$\Q$ a cyclic and separating vector for $\M$, $\N$ and $\M\cap\N$.
Consider the following properties:
\begin{itemize}
\item[$(i)$] $(\N\cap\M\subset \M,\Q)$ is -hsm,

\item[$(ii)$] $(\N'\cap\M\subset \M,\Q)$ is +hsm,

\item[$(iii)$] $J_\N J_\M =\pm J_\M J_\N \ ,$
and $J_\N \Delta_\M^{it}= \Delta_\M^{-it}J_\N$
\item[$(iii')$] $J_\N J_\M =\pm
J_\M J_\N $ and the unitaries $z(t)$ defined by
$J_\N \Delta_\M^{it}= z(2t)\Delta_\M^{-it}J_\N$
satisfy $\textnormal{Ad}z(t)\M=\M$, $\textnormal{Ad}z(t)\N=\N$,
$t\in \mathbb R$.
\end{itemize}
Property $(iii)$ holds in particular if $J_\N \M J_\N=\M$ and this,
together with properties $(i)$ and $(ii)$ characterizes \emph{local}
conformal nets \cite{Wie,AZ}. The following variation holds.
\begin{Thm} Let $\N$ and $\M$ and $\Q$ satisfy properties $(i)$, $(ii)$
and $(iii')$ above. There exists a unique conformal net $\A$ on $S^1$
such that $\A(I_1)=\M$, $\A(I_2)=\N$, with $I_1$ and $I_2$ the upper
and right semicircles, and $\Q$ is the vacuum vector. Moreover
$z_\A (t)=z(t)$.

All conformal nets on $S^1$ arise in this way.
\end{Thm}
\begin{proof} For simplicity we assume $z(t)=1$, the general case can
be treated similarly.

Set $\hat \M= J_\N \M J_\N$ so that, by property $(iii)$
$$
\Delta_{\hat\M}^{it}=\Delta_\M^{it},\quad J_{\hat\M}=\pm J_{\M}\ ,
$$
and notice that:
\begin{itemize}
\item[(a)]  $(\N\cap\M\subset \M,\Q)$ is -hsm ,
\item[(b)]  $(\N'\cap\hat \M\subset\hat\M,\Q)$ is +hsm ,
\item[(c)]  $(\N'\cap\hat \M\subset \N'\vee\M',\Q)$ is -hsm.
\end{itemize}
These are analogous to the corresponding properties in the proof of
Theorem 3 of \cite{Wie}. $(a)$ is just $(i)$, $(b)$ follows from $(a)$
by applying Ad$J_{\N}$, and $(c)$ follows by some elementary modification
of Wiesbrock proof in \cite{Wie} where, from formula $(5)$
up to formula $(8)$, one
replaces $\N'\cap \M$ with $\N'\cap\hat \M$.

It follows that the modular unitary groups $\Delta^{it_1}_{\N\cap\M}$,
$\Delta^{it_2}_{\N'\cap\hat \M}$ and
$\Delta^{it_3}_{\M}$ mutually have the same commutation relations
as the one-parameter subgroups $\Lambda_{I_1}$, $\Lambda_{I_2}$
and $\Lambda_{I_3}$ of $PSL(2,\mathbb R)$, with $I_1$ the upper-right
quarter-circle, $I_2$ the lower-right
quarter-circle, and $I_3$ the right half-circle. Therefore
the unitary group generated by these three modular unitary groups
provide a representation $U$ of the universal cover $G$ of $PSL(2,\mathbb
R)$, by an argument analogous to the one given in the proof of
\cite {GLW}, Theorem 1.2.

Again analogously to \cite{Wie}, formula $(12)$, we can see that
$U(\pi)=J_\N J_\M$,
hence
$$
U(2\pi)=J_\N J_\M J_\N J_\M =\pm 1\ ,
$$
namely $U$ is a representation of  $SL(2,\mathbb
R)$.

Set $\A(I_0) = \N$, with $I_0$ the upper semi-circle and
$\A(I)=\A(gI_0)$ if $I\in \I$ and $g\in G$ satisfy $I=gI_0$. Then
$\A(I)$ is well-defined and isotonous as in the proof of \cite{GLW},
Theorem 1.2. The rest follows by standard arguments.

That all conformal nets arise in this way follows by Theorem
\ref{genprop}.
\end{proof}
\emph{Remark.} By relaxing the condition $J_\M J_\N=\pm J_\N J_\M$
to $J_\M J_\N=\mu J_\N J_\M$ for some $\mu\in\mathbb T$, one obtains
a characterization of conformal nets on covers of $S^1$.
\section{Maximal temperature and the split property.
Examples from free probability}
Recall now that a net $\A$ is said to satisfy the {\it split
 property} if
 there exists an intermediate type I factor
 $\A(I_1)\subset F\subset\A(I_2)$ whenever the closure of the interval
 $I_1$ is contained in the interior of the interval $I_2$ \cite{DL}.

A weak form of this is the \emph{ distal split} property stating that
for each $I\in\I$ there
 exists a $\tilde I \supset I$ and a type I factor $F$ such that
 $\A(I)\subset F\subset\A(\tilde I)$.

 The following Proposition may be traced back to old argument of Kadison
 and has been used in \cite{L5,BDF}.
 \begin{Prop} If the split property holds, then the local algebras
 $\A(I)$ are approximately finite-dimensional.
 \end{Prop}
 \begin{proof}
 By continuity we may suppose that $I$ is open.
 Let $I_1\subset I_2\subset\cdots\subset I$ be an increasing sequence of open
 intervals with $\bar I_k\subset I_{k+1}$  and $\cup I_k=I$ and choose
 type I factors $\A(I_k)\subset F_k\subset \A(I_{k+1})$. Then $\A(I)$
 is generated by the increasing sequence of type I factors $F_k$,
 hence it is approximately finite-dimensional.
 \end{proof}
 Let $\A$ be a conformal net and $L$ be its conformal Hamiltonian.
 We shall say that $\A$ satisfies the {\it trace class condition} if
 $$
 \textnormal{Tr}(e^{-\b L})<\infty\ ,\forall \b > 0.
 $$
 \begin{Thm} If $\A$ satisfies the trace class condition, then
 $\A$ is split.
 \end{Thm}
 \begin{proof} By \cite{BDL}, Proposition 4.1, it is sufficient to
 show that $e^{-\b L}$ has order $0$ for all $\b>0$ and this follows
 if $e^{-\b L}$ is of type $l^p$ for all $p>0$, \cite{BDL} Lemma 2.1.
 The conclusion is thus a consequence of the following elementary lemma.
 \end{proof}
\begin{Lemma}  If $\textnormal{Tr}(e^{-\b L})<\infty\ ,\forall \b > 0$, then,
for any fixed $\b>0$,
$e^{-\b L}$ is of type $l^p$ for all $p>0$.
\end{Lemma}
\begin{proof} Let $\nu_n$ be the multiplicity of the eigenvalue
$2\pi n$ of $L$. Then the $p$-norm of $e^{-\b L}$ is
$$
|| e^{-\b L} ||_p = (\sum_n e^{-\b2\pi np}\nu_n)^{\frac{1}{p}}
=\textnormal{Tr}(e^{-\b pL})^{\frac{1}{p}}
$$
therefore the trace class condition implies $|| e^{-\b L} ||_p <\infty$
for all $p>0$.
\end{proof}
{\it Remark.} By the Kohlbecker's Tauberian theorem \cite{RV}
the trace class condition $\textnormal{Tr}(e^{-\b L})<\infty,\
\forall \b > 0$,
sets bounds for the growth of $\textnormal{Tr}(e^{-\b L})$ as $\b \to
0^+$.
\smallskip
 We now recall the following Lemma that we will need here below.
 \begin{Lemma}\label{reg} Let $V$ be a non-trivial positive energy irreducible
 unitary representation of $PSL(2,\mathbb R)$. Then
 \item[$(i)$] The restriction of $V$ to the upper triangular
 (``$ax+b$'') group is irreducible.
 \item[$(ii)$] The one-parameter unitary group
 $V\!\!\begin{pmatrix}e^t&0\\0&e^{-t}\end{pmatrix}$ is unitarily
 equivalent to the regular representation of $\mathbb R$ on
 $L^2(\mathbb R)$.
 \end{Lemma}
 \begin{proof} $(i)$: See e.g. \cite{GLW}, comments after Theorem 2.1.
 Concerning $(ii)$ we recall that, as the
logarithm of the generator of the translation unitary group
 and the generator of the dilation unitary group
satisfies the canonical commutation relation and are jointly
irreducible, the result follows by
von Neumann uniqueness theorem.
\end{proof}
 \subsection{Boltzmann statistics and maximal temperature}
 We now apply  results in the previous sections to the construction
 of a non-split net that satisfies a trace class condition with
 maximal temperature.

 Let $\H$ be a complex Hilbert space and $\F(\H)$ the Fock space over
 $\H$ with Boltzmann statistics
 $$
 \F(\H)=\bigoplus_{k=0}^{\infty}\H_k\ ,
 $$
 where $\H_{0}=\mathbb C\Q$, $\Q$ the vacuum vector, and $\H_k$
 is the $k$-fold tensor product $\H\otimes\H\cdots\otimes\H$.

 Let $\ell(h)$ and $r(k)$, $h,k\in\H$ be the left and right creation operator
 \begin{gather*}
 \ell(h)|_{\H_k}: \f_1\otimes\cdots\otimes\f_k\to
 h\otimes\f_1\otimes\cdots\otimes\f_k\\
 r(h)|_{\H_k}: \f_1\otimes\cdots\otimes\f_k\to
\f_1\otimes\cdots\otimes\f_k\otimes h
 \end{gather*}
 and
 \begin{equation*}
 	s(h)=\ell(h)+\ell(h)^*,\quad d(h)=r(h)+r(h)^*
 \end{equation*}
 the right and left fields.
 We have the following commutation relations:
 \begin{align*}
 	[\ell(h),r(k)]&=0\\
 	[\ell(h)^*,r(k)^*]&=0\\
 	[\ell(h),r(k)^*]&=(h,k)P_\Q\\
 	[\ell(h)^*,r(k)]&=(k,h)P_\Q
 \end{align*}
  with $P_\Q$ the one-dimensional projection onto $\mathbb C\Q$,
  hence
 \begin{equation}
 	[s(h),d(k)]=2i\text{Im}(h,k)P_\Q\ .
 	\label{comm}
 \end{equation}
If $H\subset\H$ is a real Hilbert space we define the von Neumann
 algebras
 $$
 \A(H)=\{s(h), h\in H\}'',\quad \B(H)=\{d(h), h\in H\}''.
 $$
 Note that $H$ is standard, i.e. $\overline{H+iH}=\H$ and
 $H\cap iH=\{0\}$, iff $\Q$ is cyclic and separating for $\A(H)$ or,
 equivalently, for $\B(H)$.

 With $Z$  the unitary involution
 $$
 Z|_{H_n}: \f_1\otimes\f_2\otimes\cdots\otimes\f_n\to
\f_n\otimes\f_{n-1}\otimes\cdots\otimes\f_1
 $$
 we have $Z\ell(h)Z=r(h)$, hence
 $$
 Z\A(H)Z=\B(H).
 $$
 \begin{Prop}[\cite{S}]\label{full} If $H$ is standard, then
 	\item[$(i)$] $\A(H)'=\B(H')=Z\A(H')Z$
 where $H'= \{k\in\H, \textnormal{Im}(h,k)=0\ \forall h\in H\}$ is the
 symplectic complement of $H$
 	\item[$(ii)$]  If the modular unitary group $\delta_H^{it}$ of $H$ on $\H$ is
 	unitarily equivalent to the regular representation of $\mathbb R$
 	on $L^2(\mathbb R)$, then $\A(H)$ is a type III$_1$-factor
 	with core isomorphic to $L(\mathbb
 	F_{\infty})\otimes B(\ell^2(\mathbb N))$.\footnote{The core a type
 	III$_1$ factor $\cal L$ is the
 	crossed product of $\cal L$ by the action of $\mathbb R$ given by the
 	modular group. $L(\mathbb F_{\infty})$ is the von Neumann algebra
 	generated by the left regular representation of the
 	free group $\mathbb F_{\infty}$ on
 	infinitely many generators. Here $\A(H)$ is
 	isomorphic to the factor considered in \cite{Rad}.}
 \end{Prop}
 We now specialize $\H$ to be the one-particle Hilbert space $\H^{(n)}$
  associated
 with the ${(n-1)}^{th}$-derivative of the $U(1)$-current algebra
  and let $H^{(n)}(I)$ be the corresponding real
 standard Hilbert subspace generated by the smooth functions with
 support in the  interval $I$ of $S^1$, see \cite{GLW}. Set
 $$
 \A_n(I)=\A(H^{(n)}(I)),\quad \B_n(I)=\B(H^{(n)}(I)).
 $$
 By the above proposition, the nets $\A_n$ are twisted local.
 \begin{Cor} $\A_n(I)$ is a type III$_1$-factor
 	with core isomorphic to $L(\mathbb
 	F_{\infty})\otimes B(\ell^2(\mathbb N))$.
 	\end{Cor}
\begin{proof} By Proposition \ref{irr} the type III$_1$ factor property
by the irreduciblity of the net which holds because $H(I)$ is standard
and $\{s(h),\ h\in\H\}''=B(\H)$.

Concerning the isomorphism class of $\A_n(I)$, by Proposition \ref{full}
it is sufficient to show that
$\delta_{H^{(n)}(I)}^{it}$ is isomorphic to the left regular representation
of $\mathbb R$. Now $\delta_{H^{(n)}(I)}^{it}$ is
unitarily equivalent to $u(\Lambda(-2\pi t))$, with $u$ the positive energy
irreducible representation of $PSL(2,\mathbb R)$ on $\H$. But this
follow by Lemma \ref{reg}.
\end{proof}
 \begin{Prop}
 	\item[(i)]  $\B_n(I)=\hat \A_n(I)$

 	\item[(ii)] $\A_n(I)\cap\B_n(I)=\mathbb C$

 	\item[(iii)] $\A_n(I)\vee\B_n(I)=B(\H)$
 \end{Prop}
 \begin{proof} $(i)$ is immediate and therefore $(ii)\:\Leftrightarrow(iii)$.
 To show $(iii)$ notice that  $\A_n(I)\vee\B_n(I)$
 contains $P_{\Q}$ by \ref{comm}, hence must be equal to $B(\H)$
 because
 $\Q$ is cyclic.
 \end{proof}
We now let $l^{(n)}$ be the the generator of the rotation one-parameter unitary
group on $\H^{(n)}$, so that the conformal Hamiltonian $L^{(n)}$ is the
promotion of
$l^{(n)}$ to $\F(\H)$, $e^{itL^{(n)}}|_{\H_k}=
e^{itl^{(n)}}\otimes e^{itl^{(n)}}\cdots\otimes e^{itl^{(n)}}$.
\begin{Prop} $\A_1$ satisfies the trace class condition with
\emph{maximal temperature} $\b_1^{-1}=
\frac{2\pi}{\log 2}$, namely
\begin{gather}
\textnormal{Tr}(e^{-\b L^{(1)}})<+\infty,\quad \b>\frac{\log 2}{2\pi}\\
\textnormal{Tr}(e^{-\b L^{(1)}})=+\infty,\quad \b\leq \frac{\log 2}{2\pi}
 \end{gather}
 \end{Prop}
 \begin{proof} Setting $L=L^{(1)}$, $l=l^{(1)}$ and $\H=\H^{(1)}$, the
occurring
 unitary representation of $PSL(2,\mathbb R)$
 on $\H$ is the irreducible one with lowest weight 1, the spectrum of
 $l$ is $\{2\pi n, n\in\mathbb N\}$ and each eigenvalue $2\pi n$ has
 multiplicity 1, hence
 $$
 \text{Tr}(e^{-\b l}) =\sum_{k=1}^{\infty}e^{-\b 2\pi k} =\frac{e^{-\b 2\pi}}
 {1-e^{-\b 2\pi}}.
 $$
 As $e^{-\b L}|_{\H_n}=e^{-\b
 l}\otimes\cdots\otimes e^{-\b l}$, we have $\text{Tr}(e^{-\b L}|_{\H_n})=
  \text{Tr}(e^{-\b l})^n$, hence we have
  $$
   \text{Tr}(e^{-\b L})=\sum_{k=0}^{\infty}\bigl(\frac{e^{-\b 2\pi}}
 {1-e^{-\b 2\pi}}\bigr)^k <+\infty \Leftrightarrow \b> \frac{\log 2}{2\pi}.
$$
 \end{proof}
 The above proposition can now  be easily generalized in order to obtain
 the following.
 \begin{Thm}\label{nonsplit} For any $T>0$, there exists a conformal
 net on $S^1$ satisfying twisted locality and the trace class condition
 with maximal temperature $\b^{-1}>T$, but not satisfying the  split
 property.
 \end{Thm}
 \begin{proof} The
 irreducible unitary representation of $PSL(2,\mathbb R)$
 on $\H^{(n)}$ is the positive energy one with lowest weight $n$ \cite{GLW}.
 We consider the associated conformal nets $\A_n$ of von Neumann algebras
on the
 full Fock space over $\H^{(n)}$. We have seen that $\A_1$ satisfies
 the trace class condition with maximal temperature $\frac{2\pi}{\log 2}$.

A similar computation with the conformal Hamitonian $L^{(n)}$ of $\A_n$
shows that $\text{Tr}(e^{-\b L^{(n)}})=\sum_{k=0}^{\infty}(\frac{e^{-\b 2\pi n}}
 {1-e^{-\b 2\pi}})^k)$ which is finite iff $\b>\b_n$ where
 $\b_n$, the solution of $\frac{e^{-\b 2\pi n}}
 {1-e^{-\b 2\pi}} =1$, satisfies $\lim_{n\to\infty}\b_n=0$.

 On the other hand $\A_n$ does not satisfies the split property  because
 by Proposition
 \ref{full} the $\A_n(I)$  have cores isomorphic to $L(\mathbb
 	F_{n})\otimes B(\ell^2(\mathbb N))$ hence they are not approximately
 	finite-dimensional.
 \end{proof}
 \emph{Remark.} It can be shown that, for each $n$, $\A_n$
 does not even satisfies the ``distal split'' property.

 \smallskip\noindent
 We conclude our paper by pointing out the following question.

 \smallskip\noindent
 {\it Problem.} Besides twisted locality, is there a notion
 related to free independence fulfilled by the nets $\A_n$?

 \smallskip
 \noindent
{\bf Acknowledgements.} We would like to thank D. Guido for
 conversations. F. R. is grateful for the hospitality
 extended to him by the CNR and the University
 of Rome ``Tor Vergata'' where this work has been done.


\begin{thebibliography}{[22]}

\bibitem{AZ} H. Araki, L. Zsido: in preparation.

\bibitem{RV} N.H. Bingham, C.M. Goldie, J.L. Teugels: ``Regular
Variations'' Cambridge Univ. Press 1987.

\bibitem{BiWi}  J. Bisognano, E. Wichmann: ``On the duality condition
for a Hermitean scalar field", J. Math. Phys. {\bf 16},  985 (1975).

\bibitem{EvBo} J. B\"ockenhauer, D. Evans: ``Modular invariants,
graphs and $\a$-inductions for nets of subfactors. I'', preprint.

\bibitem{Bo} H.-J. Borchers: ``The CPT Theorem in two-dimensional
theories of local observables", Commun. Math. Phys. {\bf 143}, 315 (1992).

\bibitem{BGL}  R. Brunetti,  D. Guido, R. Longo: `` Modular
structure and duality in conformal Quantum Field Theory'',
Commun. Math. Phys. {\bf 156}, 201-219
(1993).

\bibitem{BDF} D. Buchholz, C. D'Antoni, K. Fredenhagen: ``The
universal structure of local algebras'',
Commun. Math. Phys. {\bf 111} 123-135 (1987).

\bibitem{BDL} D. Buchholz, C. D'Antoni, R. Longo: ``Nuclear maps
and modular structures II: application to quantum field theory'',
Commun. Math. Phys. {\bf 129}, 115-138 (1990).

\bibitem{BuWi}  D. Buchholz, E. Wichmann: ``Causal independence and the
energy-level density of states in local quantum field theory'',
Commun. Math. Phys. {\bf 106}, 321-344 (1986).

\bibitem{DL} S. Doplicher, R. Longo: ``Standard and split inclusions
of von Neumann algebras'',
Invent. Math. {\bf 75} 493-536 (1984).

\bibitem{FrJo} K. Fredenhagen, M. J\"orss: ``Conformal
Haag-Kastler nets, pointlike localized fields and the existence of
operator product expansion'', Commun. Math. Phys. {\bf 176},  541 (1996).

\bibitem{FG} J. Fr\"ohlich, F. Gabbiani:
``Operator algebras and conformal field theory''
Commun. Math. Phys. {\bf 155}, 569 (1993).

\bibitem{GuLo}  D. Guido, R. Longo: ``The conformal
spin and statistics theorem'', Commun. Math. Phys.
{\bf 181},  11 (1996).

\bibitem{GLW}  D. Guido, R. Longo, H.-W. Wiesbrock: ``Extensions of
conformal nets and superselection structure'', Commun. Math. Phys.
{\bf 192}, 217-244 (1998).

\bibitem{Haag} R. Haag: ``Local Quantum Physics", Springer-Verlag, New
York-Berlin-Heidelberg 1996.

\bibitem{HiLo}  P.D. Hislop, R. Longo: ``Modular structure of the local
algebras associated with the free massless scalar field theory",
Commun. Math. Phys. {\bf 84}, 71 (1982).

\bibitem{L5} R. Longo: ``Algebraic and modular structure of von
Neumann algebras of physics'', Proc. Sympos. Pure Math. {\bf 38},
551-566 (1982).

\bibitem{Rad}  F. Radulescu: ``A one-parameter group of
automorphisms of $L(\mathbb F_{\infty})\otimes B(\H)$ scaling the
trace'', C. R. Acad. Sci. Paris {\bf 314}, no. 1, 1027-1032 (1992).

\bibitem{S} D. Shlyakhtenko: ``Free quasi-free states'', C.R. Acad. Sci.
Paris {\bf 322} S\'erie I, 645-649 (1996).

\bibitem{VDN} D.V. Voiculescu, K. Dykema, A. Nica: ``Free Random
Variables'', CRM monograph series, AMS, Providence 1992.

\bibitem{Wie} H.W. Wiesbrock: ``Conformal quantum field theory and
half-sided modular inclusions of von Neumann algebras'', Commun.
Math. Phys. {\bf 158}, 537 (1993).

\end{thebibliography}
\end{document}